\documentclass[12pt]{article}

\usepackage{amsfonts}  
\usepackage{graphicx}

\title{Eigenfunctions  for singular fully non linear equations   in unbounded regular domain } 

\author{I. Birindelli, F. Demengel}

\date{}

\catcode`@=11 \@addtoreset{equation}{section} \catcode`@=12

\newtheorem{theo}{Theorem}[section]
\newtheorem{prop}[theo]{Proposition}
\newtheorem{rema}[theo]{Remark}

\newtheorem{cor}[theo]{Corollary}
\newtheorem{lemme}[theo]{Lemma}

\def\R{{\rm I}\!{\rm  R}}
\def\grad{\nabla}

\begin{document}
\maketitle
\section{Introduction}

In this paper we  prove the existence of a generalized eigenvalue and a corresponding eigenfunction for fully nonlinear operators singular or degenerate, homogeneous of degree $1+\alpha$, $\alpha > -1$  in non bounded domains of $R^N$.  One key argument  will be the  Harnack inequality. 
 
 Very recently Davila, Felmer and Quaas \cite{DFQ, DFQ1} proved Harnack inequality in all dimensions $N$ but in the singular case  i.e. $\alpha < 0$. 
 We extend their result to the degenerate elliptic case i.e. $\alpha >0$ but only in dimension two. The proof we give uses  in an essential way this dimensional restriction. It follows the lines of the original proof of Serrin \cite{S}  in the linear case. For  Harnack inequalities in quasi-linear cases see \cite{S1} and \cite{T}. Very recently C. Imbert \cite{Im} has proved an Harnack inequality for fully-nonlinear degenerate elliptic operators; let us mention that the class of operators he considers does not include those treated in this paper (see also \cite{D} for degenerate elliptic equations in divergence form).

 It is well known that Harnack's inequality is important to control the oscillations of the solutions and hence to prove  uniform H\"older's estimates. It has been generalized to many 'weak' and nonlinear  context, we are in particular thinking of those due to Krylov and Safonov for "strong solutions" \cite{KS}, or the result of Caffarelli, Cabr\'e \cite{CC} for fully non linear equations that are uniformly elliptic. 
 Let us mention that in previous works on singular or degenerate fullynonlinear operators \cite{BD2, BD3} we proved H\"older's regularity of the solutions of Dirichlet problems in bounded domains. There the proof  relied on the regularity of the solution on the boundary and the supremum of the solution. Hence in unbounded domains that tool cannot be used.

In the case treated here of fully nonlinear operators homogenous of degree $1+\alpha$,  the Harnack inequality, due to Davila, Felmer and Quaas \cite{DFQ}, is the following 
 
  {\it Suppose that $F$ does not depend on $x$ and satisfies

  $(H1)$ and (H2) as defined later and that $-1<\alpha \leq 0$.  Suppose that $b$,  $c$ and $f$  are continuous  and that   $u$ is a nonnegative solution of 
  $$F(\nabla u, D^2 u) + b(x) \cdot \nabla u |\nabla u|^\alpha + cu^{1+\alpha} = f$$
  in $\Omega$. 
  Then  for all $\Omega^\prime \subset \subset \Omega $ there exists some constant $C$ which depends on $a$, $A$, $\alpha$, $b$, $c$, $N$, $\Omega^\prime$, $\Omega$, such that }
  $$\sup_{\Omega^\prime} u \leq C(\inf_{\Omega^\prime} u + ||f||_{L^N(\Omega^\prime)}^{1\over 1+\alpha}).$$

 Among all the  consequences of Harnack's inequality,     Berestycki, Nirenberg and Varadhan in their acclaimed paper \cite{BNV}  proved the existence of an eigenfunction for a linear, uniformly elliptic operator when no regularity of the boundary of the domain is known. The idea being that,  close to the boundary,  the solutions are controlled by the maximum principle in "small" domains, and, in the interior, one can use Harnack's inequality.

\medskip
As it is well known, inspired by \cite{BNV}, the concept of eigenvalue   in the case of bounded regular domains has lately been extended to fully-non linear operators (see \cite{BEQ}, \cite{QS}, \cite{BD2,BD3}, \cite{IY}). Two "principal eigenvalues" can be defined as the extremum of the values for which the maximum principle or respectively the minimum principle holds. 

In this article   we want to use the Harnack's inequality obtained here and in \cite{DFQ, DFQ1} to study the eigenvalue problem in unbounded domains.
Let us recall that in general, even for the Laplacian operator,  the maximum principle does not hold 
in unbounded  domain, hence we cannot  define the "principal" eigenvalue in the same way as in the case of bounded domains. In \cite{CLV} and \cite{CV} Capuzzo Dolcetta, Leoni and Vitolo study the conditions on the domain $\Omega$ in order for the Maximum principle to hold for fullynonlinear operators, extending the result of Cabr\'e \cite{Cab}.

Furthermore let us mention that in unbounded domains there are several   definitions  that allow to construct different "eigenvalues" as  the reader can  see in Berestycki and Rossi \cite{BR} for the Laplacian case. 
 Here we define the first eigenvalue as the infimum of the first eigenvalues  for bounded smooth domains included in $\Omega$. We prove the existence of  a positive eigenfunction for this so called eigenvalue, using  Harnack's inequality.

We shall also prove the existence of solutions for equations below the eigenvalues. Observe that differently from the case of bounded domain,  we can't use the maximum principle since in general it won't hold, hence again the Harnack inequality will play a key role.

\section{Assumptions on $F$}

The following hypothesis will be considered, for $\alpha>-1$:

\begin{itemize}

\item[(H1)] 
 $F$ is continuous  on $  \Omega\times\R^N\setminus\{0\}\times S\rightarrow\R$,  
and  $\forall t\in \R^\star$, $\mu\geq 0$,
 $F(x, tp,\mu X)=|t|^{\alpha}\mu F(x, p,X)$.

\item[(H2)]
For $p\in \R^N\backslash \{0\}$, $M\in S$,  $N\in S$, 
$N\geq 0$

\begin{equation}\label{eqaA}
a|p|^\alpha tr(N)\leq F(x,p,M+N)-F(x,p,M) \leq A
|p|^\alpha tr(N).
\end{equation}

\item[(H3)]
There exists a continuous function $\tilde \omega$, $\tilde \omega(0)=0$
such that for all $(x,y)\in \Omega^2 $,  $\forall p\neq 0$, $\forall X\in S$
$$|F(x,p,X)-F(y,p,X)|\leq \tilde \omega(|x-y|) |p|^\alpha |X|.$$
\item [(H4)]
 There exists a
continuous function $  \omega$ with $\omega (0) = 0$, such that if
$(X,Y)\in S^2$ and 
$\zeta\in \R^+$ satisfy
$$-\zeta \left(\begin{array}{cc} I&0\\
0&I
\end{array}
\right)\leq \left(\begin{array}{cc}
X&0\\
0&Y
\end{array}\right)\leq 4\zeta \left( \begin{array}{cc}
I&-I\\
-I&I\end{array}\right)$$
and $I$ is the identity matrix in $\R^N$,
then for all  $(x,y)\in \R^N$, $x\neq y$
$$F(x, \zeta(x-y), X)-F(y,  \zeta(x-y), -Y)\leq \omega
(\zeta|x-y|^2).$$

\end{itemize}

Observe that when $F$ is independent of $x$,  conditions (H3) and (H4) are not needed. 

\begin{rema} When no ambiguity arises we shall sometimes write 
$F[u]$ to signify $F(x,\grad u,D^2u)$.
\end{rema}

Recall that examples of operators satisfying these conditions include the $p$-Laplacian with $\alpha=p-2$ and
$$F(\grad u,D^2u)=|\grad u|^\alpha{\cal M}_{a,A}^\pm (D^2u)$$
where ${\cal M}_{a,A}^+$ is the Pucci operator
${\cal M}_{a,A}^+(M)=A{\rm Tr}(M^+)-a{\rm Tr}(M^-)$.
and ${\cal M}_{a,A}^-(M)=a{\rm Tr}(M^+)-A{\rm Tr}(M^-)$.

We assume that $h$  and $V$ are some continuous bounded functions on $\bar\Omega$ and 

\noindent (H5)
- Either $\alpha \leq 0$ and $h$  is H\"older continuous  of exponent $1+\alpha$,

- or  $\alpha >0$ and 
$$[(h(x)-h(y))\cdot (x-y)]\leq 0$$

The solutions that we consider will be taken in the sense of viscosity, see e.g. \cite{BD1} for precise definitions, let us recall that in particular we do not test when the gradient of the test function is null .

\section{Main results}
\subsection{The Harnack's inequality in the two dimensional case.}
In this subsection we state the Harnack's inequalities that will be proved in section 5 and used in section 4, together with some important corollary.

\begin{theo}[Harnack's inequality] \label{prop1}
Suppose that $\Omega$ is  a bounded domain in $\R^2$, and that $F$  satisfies (H1) to (H4), $h$ satisfies (H5). 

Let $u$ be a positive solution of
\begin{equation}\label{eqhar}
F(x,\nabla u, D^2 u) +h(x).\nabla u|\nabla u|^\alpha+ V(x) u^{1+\alpha}= 0  \ \mbox{ in}\ \Omega .
\end{equation}
Let $\Omega^\prime \subset\subset \Omega $. Then there exists $K=K(\Omega,\Omega^\prime, A,a, |h|_\infty, |V|_\infty  )$ such that
\begin{equation}\label{esthar}\sup_{\Omega^\prime}  u\leq K\inf_{\Omega^\prime} u .
\end{equation}
\end{theo}

\begin{theo}[Harnack's inequality] \label{prop2}
Under the same hypothesis of Theorem \ref{prop1}, for  $f$ a  bounded continuous function on $\Omega$, let $u$ be a positive solution of
\begin{equation}\label{eqhar2}
F(x,\nabla u, D^2 u) +h(x).\nabla u|\nabla u|^\alpha=  f(x) \ \mbox{ in}\ \Omega .
\end{equation}
Let $\Omega^\prime \subset\subset \Omega $. Then there exists $K=K(\Omega,\Omega^\prime,A,a, |h|_\infty )$ such that
\begin{equation}\label{esthar2}\sup_{\Omega^\prime}  u\leq K\left(\inf_{\Omega^\prime} u +|f|_{L^\infty(\Omega)} ^{1\over 1+\alpha}\right).
\end{equation}
\end{theo}
\begin{rema}\label{remV}
The result in theorem \ref{prop2}  still holds for $u$ a positive solution of
$$
F(x,\nabla u, D^2 u) +h(x).\nabla u|\nabla u|^\alpha+V(x)u^{1+\alpha}=  f(x) \ \mbox{ in}\ \Omega .
$$with 
$V$ continuous,  bounded  and $V\leq 0$. In that case the constant $K$ depends also on $|V|_\infty$. 
\end{rema}

\begin{cor}\label{corhold}
Let $u$ be a solution of (\ref{eqhar2}). 
Let $R_o$  be such that $B(0, R_o)\subset \Omega$. Then there exists  $K$  which depend only on  $ A,a, |h|_\infty$ and $R_o$,  such that  for any $R<R_o$:
\begin{equation}\label{estharR}\sup_{B(0, R)}  u\leq K(\inf_{B(0, R)} u+  R^{2+\alpha\over 1+\alpha}|f|_{L^\infty(B(0, R_o)} ^{1\over 1+\alpha}).
\end{equation}

As a consequence, for any  solution $u$  of (\ref{eqhar2}), for any $\Omega^\prime\subset\subset\Omega$ there exists $\beta\in (0,1)$ depending on the Harnack's constant in (\ref{estharR}) such that  
$u\in C^{o,\beta}(\Omega^\prime)$.
\end{cor}

 \medskip

An immediate consequence of Harnack's inequality is the following Liouville type result  : 

\begin{cor}[Liouville]\label{Liou}
 Let $u$ be a solution of 
$F(x,\nabla u, D^2 u)= 0$ in $\R^2$,  if $u$ is bounded from below, then $u\equiv cte$.
\end{cor}
See \cite{CL} for other Liouville results.

\subsection{Existence's  results in unbounded domains.}
Before stating the results in unbounded domains we  recall what we mean by first eigenvalue and the property of these eigenvalues in the  bounded case. 

When $\Omega$ is a bounded domain we define 
$$\overline\lambda (\Omega)= \sup \{ \lambda, \exists \ \varphi>0\ {\rm in} \ \Omega , F[\varphi]+h(x)\cdot\nabla \varphi|\nabla \varphi|^\alpha+  (V(x)+ \lambda) \varphi^{1+\alpha} \leq 0\}$$
and
$$\underline\lambda (\Omega)= \sup \{ \lambda, \exists \ \varphi<0\ {\rm in} \ \Omega , F[\varphi]+h(x)\cdot\nabla \varphi|\nabla \varphi|^\alpha+  (V(x)+ \lambda) \varphi|\varphi|^{\alpha} \geq 0\}.$$
When $\Omega$ is a {\bf bounded} regular domain, we proved in \cite{BD1} that 
there exists $\varphi>0$  and $\psi<0$ in $\Omega$ which are  respectively a solution of 
$$\left\{ \begin{array}{cc}
F(x,  \nabla \varphi(x),
D^2\varphi(x))+h(x)\cdot\nabla \varphi|\nabla \varphi|^\alpha+  (V(x)+\overline\lambda (\Omega)) \varphi^{1+\alpha}=0 & {\rm in }\ \Omega\\
 \varphi = 0  & {\rm on}\ \partial \Omega
 \end{array}\right.$$
 and
 $$\left\{ \begin{array}{cc}
F(x, \nabla \psi(x),
D^2\psi(x))+h(x)\cdot\nabla \psi|\nabla \psi|^\alpha+  (V(x)+\underline\lambda (\Omega))|\psi|^{\alpha}\psi=0 & {\rm in }\ \Omega\\
 \psi= 0  & {\rm on}\ \partial \Omega.
 \end{array}\right. $$

Moreover $\varphi$ and $\psi$ are  H\"older continuous.

 We assume that $\Omega\subset\R^N$ is  a ${\cal C}^2$  possibly unbounded domain. Define
 
$$\overline\lambda (\Omega)=\inf\{\overline\lambda(A),\mbox{for all smooth bounded domain } A, A\subset \Omega\},$$
and
$$\underline\lambda (\Omega)=\inf\{\underline\lambda(A),\mbox{for all smooth bounded domain } A, A\subset \Omega\}.$$

When no ambiguity arises we shall omit to write the dependence of the eigenvalues with respect to the set $\Omega$.

We wish  first  to give some bounds on  $\overline\lambda(\Omega)$.
For simplicity this will be done for  $h\equiv 0$, $V\equiv 0$. 
If $\Omega$ is bounded it is easy to see that $\overline\lambda(\Omega)>0$, while it is obvious that for $\Omega=\R^N$, $\overline\lambda(\Omega)=0$. 
We wish to prove that this is not the case for all unbounded domains, in fact we shall see that, as long as $\Omega$ is bounded in one direction, then  $\overline\lambda(\Omega)>0$.

\begin{prop} Suppose that $\Omega$ is contained in a strip of width $M$ i.e. up to translation and rotation
$$\Omega\subset [0,M]\times\R^{N-1}$$
then there exists $C=C(\alpha,a,A)>0$ such that
\begin{equation}\label{estiml}
\overline\lambda(\Omega)\geq\frac{C}{M^{2+\alpha}}.
\end{equation}
\end{prop}
{\bf Proof:} 
Fixe $\gamma\in (0,1)$ and observe that 
$u(x)=\sin^\gamma(x_1\frac{\pi}{4M}+\frac{\pi}{8})\geq 0$ in $\Omega$ and
\begin{eqnarray*}
F[u]
&\leq&\gamma^{\alpha+1}\left(\frac{\pi}{4M}\right)^{\alpha+2} \sin^{\gamma(\alpha+1)-(\alpha+2)}(x_1\frac{\pi}{4M}+\frac{\pi}{8})\left(\cos(x_1\frac{\pi}{4M})+\frac{\pi}{8}\right)^\alpha\cdot\\
&& a [\gamma-1-\gamma\sin^2(x_1\frac{\pi}{4M}+\frac{\pi}{8})].
\end{eqnarray*}

Hence, using  
$$\frac{\pi}{2}>\frac{3\pi}{8}\geq x_1\frac{\pi}{4M}+\frac{\pi}{8}>\frac{\pi}{8}>0,$$ 
 we get that there exists $C=C(\gamma,a,\alpha)$

$$F[u]+\frac{C}{M^{2+\alpha}} u^{\alpha+1}\leq 0.$$
Clearly this implies that 
$\overline\lambda(A)\geq\frac{C}{M^{2+\alpha}} $ for any $A\subset \Omega$. This gives (\ref{estiml}) and it ends the proof.

In the next theorem we want to be in the same hypothesis for which   Harnack's  inequality holds, hence we consider the following condition:

\noindent(C) {\em $F$ satisfies (H1), (H2);
 if $N\geq 3$ $F$ is independent of $x$ and $-1<\alpha\leq 0$; if $N=2$, $\alpha>-1$,  $F$ may depend on $x$ and it satisfies (H3) and (H4) .}

 \begin{theo}
 \label{exieig}
 Suppose that $\Omega$ is some smooth domain possibly non bounded, of $\R^N$.  Suppose that $F$ satisfies (C), that $h$ satisfies (H5), and that $V$ is continuous,  and bounded.
 Then there exist some  functions $\phi>0$ and  $\psi<0$  which are continuous and satisfy, respectively
 $$F[\phi]+ h(x)\cdot \nabla \phi |\nabla \phi|^\alpha + (\overline\lambda (\Omega) + V(x)) \phi^{1+\alpha}  = 0\ {\rm in} \ \Omega,$$
 $$F[\psi]+ h(x)\cdot \nabla \psi |\nabla \psi|^\alpha + (\underline\lambda (\Omega) + V(x)) |\psi|^{\alpha} \psi = 0\ {\rm in} \ \Omega.$$
 Furthermore $\phi$ and $\psi$ are H\"older continuous.
  \end{theo}
  
  In the next proposition we treat existence of solutions below the eigenvalues.
\begin {prop}\label{below}
For any $\lambda<\underline\lambda(\Omega)$, for any 
$f\in {\cal C}_c (\Omega) $ non positive,  there exists $v> 0$ solution of

 $$F[ v]+h(x)\cdot \nabla v|\nabla v|^\alpha +(\lambda +V(x))v^{1+\alpha}=f\  \mbox{in}\ \Omega.$$
 
Furthermore, for $f\not \equiv 0$ there exists $C$, which depends on the support of $f$, such that
$$\vert v\vert _\infty\leq C\vert f\vert _\infty^{\frac{1}{1+\alpha}}.$$
Similarly if  $\lambda<\underline\lambda(\Omega)$, for any 
$f\in {\cal C}_c (\Omega)\geq 0$, there exists $v< 0$ solution of

 $$F[ v]+h(x)\cdot \nabla v|\nabla v|^\alpha +(\lambda +V(x))|v|^{\alpha}v=f\  \mbox{in}\ \Omega.$$

\end{prop}

\begin{rema}
As mentioned in the introduction, in \cite{BD2} we proved some H\"older's regularity result  for all $\beta \in [0, 1[$  in bounded regular domains,  see Proposition \ref{prophold}, but for homogeneous or regular boundary conditions. More precisely  the H\"older's constants depend on the $L^\infty$ norm  of $u$ and $u$ is zero on the boundary. From this we derive some H\"older's  uniform estimates for sequences of solutions and this allows to prove that a sequence of  such solutions converges for a subsequence  towards a solution. 
This cannot be used in the proof of the results above, indeed  we shall need compactness results inside bounded sets  $\Omega_n$ whose size increases, for sequence of functions which have  uniform $L^\infty$ bounds   on bounded fixed sets, but for which $L^\infty(\Omega_n) $ norm may  go to infinity.  
\end{rema}

%We now recall the  properties of the first eigenfunctions : 

%\begin{theo}\label{pmaxreg}
 % Suppose that $\lambda < \overline\lambda (\Omega) $ and that $u$ satisfies 
  %$$\left\{ \begin{array}{cc}
  %F(x, \nabla u, D^2 u) + h(x) \cdot \nabla u  |\nabla u|^\alpha + (V(x)+\lambda) |u|^\alpha u \geq 0&{\rm in } \ \Omega\\
  %u\leq 0 \ {\rm on } \ \partial \Omega
  %\end{array}\right.$$
  %Then $u\leq 0$ in $\Omega$.
  %\end{theo}
  
 \section{Known results.}
   
 We now recall the following weak  comparison principle  which will be used for the proof of Theorem \ref{prop2}.
 
 \begin{theo}\label{thcomp1}
 Suppose that  $F$,  $h$ and $V$ are as above and that $V\leq 0$. 

Suppose that $f$ and $g$ are continuous and bounded and that 
$u$ and
$v$  satisfy
\begin{eqnarray*}
 F(x, \nabla u, D^2 u)+ h(x)\cdot \nabla u |\nabla u|^\alpha +V(x)|u|^\alpha u& \geq & g\quad 
\mbox{in}\quad \Omega \\ 
F(x,  \grad v,D^2 v)
+ h(x)\cdot \nabla v |\nabla v|^\alpha+ V(x) |v|^{\alpha}v & \leq & f \quad  \mbox{in}\quad
\Omega  \\ 
u \leq  v &&   \quad  \mbox{on}\quad \partial\Omega.
\end{eqnarray*}
 Suppose that $f< g$, then $u \leq v$ in $\Omega$. 
Moreover if $V<0$ and $f\leq g$ the result still holds. 

\end{theo}
We shall also need for the proof of Theorem  \ref{prop1} another comparison principle : 
\begin{theo}\label{thcomp2}
Suppose that $\tau< \overline\lambda (\Omega)$, $f\leq 0$, $f$ is
upper semi-continuous and
$g$ is lower semi-continuous  with $f\leq  g$.

Suppose that there exist
$u$  continuous and
$v\geq 0$ and   continuous, satisfying
\begin{eqnarray*}
 F(x, \nabla u, D^2 u)+ h(x)\cdot \nabla u |\nabla u|^\alpha +(V(x)+\tau) |u|^\alpha u& \geq & g\quad 
\mbox{in}\quad \Omega \\ 
F(x,  \grad v,D^2 v)
+ h(x)\cdot \nabla v |\nabla v|^\alpha+( V(x)+\tau) v^{1+\alpha} & \leq & f \quad  \mbox{in}\quad
\Omega  \\ 
u \leq  v &&   \quad  \mbox{on}\quad \partial\Omega.
\end{eqnarray*}
Then $u \leq v$ in $\Omega$ in each of these two cases:

\noindent 1) If $v>0$ on $\overline{ \Omega}$ and either $f<0$ in 
$\Omega$, 
 or  $g(\bar x)>0$ on every point $\bar x$ such that  $f(\bar x)=0$.

\noindent  2) If $v>0$ in $\Omega$, $f<0$ in $\overline{\Omega}$  and $f<g $ on
$\overline\Omega$. 
\end{theo}

The proof can be found in \cite{BD1}. 
We also recall some  regularity results 
\begin{prop}\label{prophold}

Suppose that $F$ satisfies (H1),(H2), (H3). 
Let $f$ be some continuous  function in $\overline{\Omega}$.
 Let $u$ be a viscosity non-negative bounded solution of 
\begin{equation}\label{eq4.1}
\left\{
\begin{array}{lc}
F(x, \nabla u, D^2u)+ h(x)\cdot\nabla u|\nabla u|^{\alpha}=f & \ {\rm in}\
\Omega\\
u=0 &  \ {\rm on}\ \partial\Omega.
\end{array}
\right.
\end{equation}
Then, if $h$ is continuous and bounded,  for
any
$\gamma<1$ there exists some constant
$C$ which depends only on $|f|_\infty$, $|h|_\infty$ and $|u|_\infty$ such that :

$$|u(x)-u(y)|\leq C|x-y|^\gamma$$
for any
$(x,y)\in\overline\Omega^2$.
\end{prop}
%This result has the following compacity consequence 

%\begin{cor} \label{corcomp}
% Suppose that $(f_n)$ is a  
%sequence of continuous and uniformly bounded functions, and $(u_n)$ is a
%sequence of bounded viscosity solutions of  

%$$\left\{\begin{array}{lc}
%F(x, \nabla u_n, D^2u_n)+ h(x).\nabla u_n|\nabla
%u_n|^{\alpha}=f_n &\mbox{in}\quad\Omega\\
%u_n=0 &\mbox{in}\quad\partial\Omega.
%\end{array}
%\right. .$$ 
%Then the sequence  $(u_n)$
%is relatively
%compact in 
%${\mathcal C} (\overline\Omega)$.
%\end{cor}

Under slightly stronger condition on $F$ we also prove the Lipschitz regularity of the solutions.

\section{Proofs of the  Main results}
\subsection{Proof of existence}

We start by proving the existence results:

 {\em Proof of Theorem \ref{exieig}.}
 We shall only explicitly write the proof of the existence of $\phi>0$, the case of $\psi<0$ being analogous.
 Let $\Omega_n$ be a sequence of bounded  subsets such that 
 $$\Omega_n\subset \subset \Omega_{n+1}\subset \subset \Omega,\quad \overline\lambda (\Omega_n)\rightarrow \overline\lambda (\Omega)\quad \mbox{ and} \quad \cup_{n}\Omega_n=\Omega.$$ Let $f_n$ be a sequence of functions in ${\cal C}_c (\Omega_n\setminus\overline{\Omega_{n-1}})$, $f_n\leq 0$ and not identically zero. 
 Since $\overline\lambda(\Omega_n)> \overline\lambda (\Omega)$,  for all $n$ there exists $u_n\geq 0$ which solves 
 $$ \left\{\begin{array}{cc}
F[ u_n]+ h(x)\cdot \nabla u_n |\nabla u_n|^\alpha + (\overline\lambda(\Omega) + V(x)) u_n^{1+\alpha} = f_n&\ {\rm in} \ \Omega_n\\
 u_n = 0&\ {\rm on }\  \partial \Omega_n.
 \end{array}\right.
 $$
 Let $x_0\in \Omega_1$, then $u_n(x_0)>0$ for all $n$ by the strict maximum principle.  Define
  
 $$v_n(x) = {u_n(x)\over u_n(x_0)}$$
 that we extend by zero outside $\Omega_n$,  obtaining  in such a way a continuous function. 
Let $O$ be a bounded regular domain in $\Omega$.  We prove that $v_n$ converges uniformly on $K= \overline{O}$ . Indeed there exists $N_0$ such that $\Omega_{n} $ contains $K$ for all $n\geq N_0$.
 As a consequence on $O$ , for $n\geq N_0$ 
  $$\begin{array}{cc}
 F[v_n]+h(x)\cdot \nabla v_n |\nabla v_n|^\alpha + ( \overline\lambda (\Omega)+V(x)) v_n^{1+\alpha} = 0&\ {\rm in} \ O.\\
 \end{array}
 $$
 Moreover $v_n(x_0)=1$.  Using  Harnack's inequality   of Theorem \ref{prop1} we know that there exists some constant $C_K$ such that
 $$ \sup v_n \leq C_K(\inf v_n)\leq C_K.$$
This implies in particular that $v_n$ is bounded independently of $n$ in $K$. 

By taking $ f_n = -V(x) v_n^{1+\alpha}$ in Corollary \ref{corhold}   on the  open set $O$,  one gets that $(v_n)_n$ is relatively compact in $O$. A subsequence of $v_n$ will  converge to a solution $\phi$ of 
 $$F[\phi]+ h(x)\cdot \nabla \phi |\nabla \phi|^\alpha + ( \overline\lambda (\Omega)+ V(x)) \phi^{1+\alpha}= 0\quad\mbox{in }\quad O.$$
  $\phi(x_0)= \lim v_n(x_0)=1$ implies that $\phi$ cannot be identically zero. 
By strict maximum principle on compacts sets  of $\Omega$, $\phi>0$ inside $\Omega$. This ends the proof.

\medskip

\medskip
{\em Proof of Proposition \ref{below}.}  We consider only the case $f\leq 0$ and $\lambda < \overline{\lambda} (\Omega)$. 
We first treat the case $f\not\equiv 0$. Let $K$ be the compact support of $f\leq 0$.
As in the previous proof let 
$\Omega_n$ be a sequence of bounded  sets such that

$$\Omega_n\subset\Omega_{n+1}\ \mbox{and}  \ \cup_n\Omega_n=\Omega.$$
Let $u_n$ be a (positive ) solution of
$$
\left\{\begin{array}{lc}
F[u_n]+h(x)\cdot \nabla u_n|\nabla u_n|^\alpha + (V(x)+\lambda) u_n^{1+\alpha}=f  & \mbox{in}\ \Omega_n\\
u_n=0 & \mbox{on}\  \partial\Omega_n.
\end{array}
\right.
$$
Let $\varphi^+$  be given in Theorem \ref{exieig} such that 
$$F[\varphi^+]+h(x)\cdot \nabla \varphi^+|\nabla \varphi^+|^\alpha + ( \overline\lambda (\Omega)+ V(x)) {\varphi^+}^{1+\alpha} = 0$$
with $L^\infty$ norm $1$ in $K$.

Rescaling $\varphi^+$,
$$\varphi_1 = {\varphi^+ \sup |f|^{1\over 1+\alpha}\over (\overline\lambda -\lambda)^{1\over 1+\alpha} \inf_K \varphi^+},$$
by homogeneity is a solution of

$$F[\varphi_1] +h(x)\cdot \nabla \varphi_1|\nabla \varphi_1|^\alpha + ( \lambda + V(x)) \varphi_1^{1+\alpha} = (\lambda-\overline\lambda) {(\varphi^+)^{1+\alpha} \sup (-f)\over (\overline\lambda-\lambda )(\inf_K \varphi^+)^{1+\alpha}} \leq  f.$$
We can apply the comparison principle  Theorem \ref{thcomp2} in $\Omega_n$, since $\varphi_1>0$ on $\partial \Omega_n$,  to derive that

$$ 0\leq u_n\leq \varphi_1$$
for any $n$.  Using the same argument as  in the proof of Theorem \ref{exieig}, on every compact subset of $\Omega$ there is a subsequence of $(u_n)_n$ converging to $u$, a  solution of
$$
F[u]+h(x)\cdot \nabla u|\nabla u|^\alpha + (V(x)+\lambda) u^{1+\alpha}=f \quad \mbox{in}\ \Omega.
$$
By the strict maximum principle  applied on bounded sets of $\Omega$ we get that $u>0$.

We now prove the case  $f\equiv 0$. Without loss of generality 
we only treat the case $\lambda<\overline\lambda(\Omega)$.

Let $\Omega_n$ be a sequence of bounded  sets such that

$$\Omega_n\subset\Omega_{n+1}\ \mbox{and}  \ \cup_n\Omega_n=\Omega.$$

Let $u_n$ be a  solution of
$$
\left\{\begin{array}{lc}
F[u_n]+h(x)\cdot \nabla u_n|\nabla u_n|^\alpha + (V(x)+\lambda) u_n|u_n|^\alpha=0  & \mbox{in}\ \Omega_n\\
u_n=1 & \mbox{on}\  \partial\Omega_n.
\end{array}
\right.
$$

Since $\lambda<\inf \{\overline{\lambda}(\Omega_n)\}$,  $u_n$ exists, is well defined
and $u_n> 0$ in $\Omega_n$. 
Let $P_0\in \Omega_1$ .

Rescaling $u_n$ we get that $v_n = {u_n\over u_n(P_0)}$ is a solution of
$$F[v_n]+h(x)\cdot \nabla v_n|\nabla v_n|^\alpha + (V(x) + \lambda) v_n^{1+\alpha} = 0.$$
By  Harnack's inequality,  for every  relatively compact domain $O$, 
$v_n$ is bounded  on $K= \overline{O}$ . 

Using the compactness results on $O$  there exists  a subsequence $v_n$ which converges uniformly to some $v$ solution of 
$$F[v]+h(x)\cdot \nabla v|\nabla v|^\alpha + (V(x) + \lambda) v^{1+\alpha} = 0.$$

Moreover,  since $v_n(P_0)=1$,  and the convergence is uniform one gets that $v(P_0)=1$, hence $v$ is not identically zero and by the strict maximum principle $v>0$ in $\Omega$.

\subsection{Proofs of Harnack's inequality in the two dimensional case.  }

The proofs that we propose follow the lines  in  Gilbarg Trudinger \cite{GT} and  Serrin \cite{S}, with some new arguments that make explicite use of the eigenfunction in bounded domains.
This  extends  the result of \cite{DFQ} to the case $\alpha >0$, but only in the two dimensional case.

In the proof of Theorem \ref{prop1}  and \ref{prop2} we shall use the following  lemma

   \begin{lemme}\label{lemme1}
 
  Suppose that  $F$, $h$ and $V$ are as above. Let $b$ and $c$,  be some positive parameters, $x_o=(x_{o1},x_{o2})\in \R^2$.  
  Let $$E= \{ x=(x_1,x_2),\  \sigma^2(x):=  {(x_1-x_{o1})^2\over b^2}+ {(x_2-x_{o2})^2 \over c^2} \leq 1, \ x_1-x_{o1}> {b\over 2}\}.$$
  Then there exists a constant $\gamma>0$ such that  
  $$v = {e^{-\gamma \sigma^2}-e^{-\gamma}\over e^{-\gamma/4}-e^{-\gamma}},$$

 satisfies in $E$
\begin{equation}\label{lm1}
F(x, \nabla v, D^2 v)-|  h|_\infty  |\nabla v|^{1+\alpha} -| V|_\infty  v^{1+\alpha} > 0.
\end{equation}
 (Note that $v$  is strictly positive inside $E$ and is zero on the elliptic part of the boundary).
 
 \end{lemme}
 
 \begin{rema}\label{rem7}
  The same  result holds for the symmetric part of ellipsis : 
   $E= \{ x=(x_1,x_2),\  \sigma^2(x)\leq 1, \ x_1-x_{o1}< {-b\over 2}\}$. 
   \end{rema}
   
Proof of Lemma \ref{lemme1}.

Without loss of generality one can assume that $x_o=0$. 

 Let $\tilde v =  {e^{-\gamma \sigma^2}\over e^{-\gamma/4}-e^{-\gamma}}$ and let $B$ be the diagonal $2\times 2$ matrix such that $B_{11}=\frac{1}{b^2}$ and $B_{22}=\frac{1}{c^2}$.
Then $\nabla v = -2\gamma Bx \tilde v$ and 
$$D^2 v = (2\gamma) (2\gamma B x\otimes B x-B) \tilde v.$$
Since $B$ and $Bx\otimes Bx$ are both nonnegative, 
 $$a(tr(D^2 v)^+)-A (tr(D^2 v)^-)\geq \left(a \gamma^2  4({x_1^2\over b^4}+ {x_2^2\over  c^4})-2(A+a)\gamma ({1\over b^2}+ {1\over c^2})\right) \tilde v.$$
We define 
$$m = \inf \left(b^{-\alpha}, 2^\alpha ({1\over b^2}+ {1\over c^2})^{\alpha/2}\right)$$
and $M = 2^{1+\alpha} ({1\over b^2} + {1\over c^2})^{1+\alpha\over 2}$.
We choose 
\begin{equation}\label{gam} \gamma =\sup \left({4(A+a)\over a} (1+  {b^2\over  c^2}), {4|h|_\infty Mb^2\over ma}, \left({4|V|_\infty b^2 \over am}\right)^{1\over 2+\alpha}\right).
\end{equation}
Using (H1):
\begin{eqnarray*}
&& F(x, \nabla v, D^2 v)+  h(x)\cdot \nabla v |\nabla v|^\alpha + V(x) v^{1+\alpha}\geq\\
&\geq&|\nabla v|^{\alpha}(a(tr(D^2 v)^+)-A (tr(D^2 v)^-)) -|h|_\infty |\nabla v|^{1+\alpha} -|V|_\infty v^{1+\alpha} > 0.
\end{eqnarray*}

This ends the proof of Lemma \ref{lemme1}.
\begin{rema}\label{rem2}

The proof in the case  $f\not\equiv 0$ follows the lines of the case $f \equiv0$ but the ellipsis are rescaled.   Hence we shall use, for $\rho_o$ to be defined, $\sigma ({x\over \rho_o})$ instead of $\sigma$. It  will be important to observe that $\gamma$ does not depend on bounded $\rho_o$.  This is  immediate from the definition of $\gamma$  in (\ref{gam}) and the constants $m$, $M$, $b$ and $c$ involved. 
\end{rema}

\noindent{\em Proof  of Theorem \ref{prop1}:}

Let us remark that the existence  of a positive solution $u$ implies in particular that $ \overline\lambda (\Omega)\geq 0$.
Moreover  without loss of generality we can suppose that $\overline\lambda(\Omega)>0$. Indeed, by the properties of the eigenvalue there exists $\Omega_1\subset \Omega $ such that $\Omega^\prime\subset\subset \Omega_1$ and $\overline\lambda(\Omega_1)>\lambda(\Omega)\geq 0$. Then we consider the proof in $\Omega_1$ instead of $\Omega$.

We shall prove the following claims : 

\noindent{\bf Claim 1:} Suppose that $\Omega= B(0,1)$.    For any $P\in B(0,{\frac{1}{3}})$ there exists $ K$ which depends only on 
$a$, $A$,  and   bounds on $h$ and $V$ such that
$$u(P)\geq Ku(0).$$
{\bf Claim 2:}  For any $P\in B_{\frac{1}{4}}(0)$,  there exist $K_1$ and $K_2$ such that
$$ K_1 u(0)\leq u(P)\leq K_2 u(0).$$
{\bf Claim 3:} Suppose that $\Omega=B(0,R)$. For any $P\in B(0,\frac{R}{4})$  such that 
$$ K_1 u(0)\leq u(P)\leq K_2 u(0),$$
where $K_1$ and $K_2$ depend on $R$ only when $h$ and $V$ are not identically 0.

\noindent{\bf Claim 4:}  The inequality holds true  for $\Omega$ bounded and $\Omega^\prime \subset\subset \Omega$. 

\medskip

Proof of Claim 1 :

So we are in the case
$\Omega = B(0,1)$ with $\overline\lambda(B(0,1))>0$. Hence there exists $\delta>0$ sufficiently small such that $\overline\lambda(B(0,1+\delta))>0$ as well.

% For some for $\delta >0$ small enough $0< \bar\lambda (B(0, 1+ \delta)) $. 
 %Indeed let $0<\epsilon < \bar\lambda (B(0,1))/2$ then there exists $v>0$ a solution of 
%$$F(\nabla v, D^2 v) + h \nabla v|\nabla v|^\alpha+ (V(x)+ \epsilon)v^{1+\alpha} = -1$$
 %Then there exists $\delta >0$ such that $v+\delta_v $ is a solution of 
 %$$F(\nabla v, D^2 v) + h \nabla v|\nabla v|^\alpha+ (V(x)+ \epsilon)(v+\delta_v)^{1+\alpha} = -1/2$$
 %The function $v+\delta_v$ is positive   on some ball $B(0,1+\delta)$ which implies that %$\bar\lambda(B(0, 1+\delta))>0$ .

Let  $u_{\delta}$ be the  corresponding positive eigenfunction  such that $u_{\delta}$ has the $L^\infty $ norm equals to  ${1\over 2}$,   i.e. $u_\delta$ satisfies   
$$\left\{ \begin{array}{lc}F[u_\delta]+  h(x)\cdot \nabla u_{\delta}|\nabla u_{\delta}|^\alpha  + \left(V(x)+ \overline\lambda(B(0,  1+\delta)) \right)u_{\delta}^{1+\alpha} = 0\ {\rm in } \ B(0,1+\delta)&\ \\
 u_{\delta}= 0\ {\rm on }\quad  \partial B(0,1+ {\delta}).&\ 
 \end{array}\right.
 $$
Let $\chi = u(0) u_{\delta}$.

Let 
 $G_1= \{ x\in B(0,1), u(x) > \chi (x)\}$. 
  The connected component of $G_1$, denoted $G$,   which contains $0$,  contains at least one point on $\partial B(0,1)$. Indeed, if not, on  the boundary of $G$ one would have $u(x)\leq \chi$ and since $0 < \overline\lambda (B(0,1+ \delta))$,  $\chi$ is a supersolution  of 
  $F[ \chi]+h(x).\nabla \chi |\nabla \chi|^\alpha+ ( V(x) ) \chi ^{1+\alpha}< 0 $,
  then applying the comparison Theorem  \ref{thcomp2}  in the set $G$, one would get 
  $u(x)\leq \chi$ inside $G$,   but   this does not hold at the point $0$  since $\sup u_{\delta}= {1\over 2}$,    so we have reached a contradiction.   Without loss of generality we will suppose that the  boundary point  has coordinates $(0,1)$.      We denote by $\kappa $ the positive constant $\displaystyle\inf_{B(0,1)} u_{\delta}$, and $\chi_1 = \kappa u(0)$

  We now  introduce the part of ellipsis  $E_i$ $i=1,2,3$  given by: 
  
  $$E_1 = \{ (x_1,x_2),\ {(x_1+{5\over 2})^2\over 9} + 4(x_2-{\sqrt{3}\over 4})^2 \leq 1, \ x_1\geq -1\}$$
  $$E_2 = \{ (x_1,x_2),\  {(x_1-{5\over 2})^2\over 9} + 4(x_2-{\sqrt{3}\over 4})^2 \leq 1, \ x_1\leq 1\}.$$
 Observe that the segment $[-1/2, 1/2] \times\{ {\sqrt{3}\over 4}\}$ is contained in $E_1\cap E_2$. 
 while  $(0,0)\not \in E_1\cap E_2\subset B(0,1)$.
 
 The third part of  ellipse  $E_3$ has its  straight part in $E_1\cap E_2$ and vertex at $(0,-1)$:

 $$E_3 = \{ (x_1,x_2), \ 4x_1^2+ \left({x_2-1-{\sqrt{3}\over 2}\over 2+ {\sqrt{3}\over 2}}\right)^2 \leq 1, x_2\leq \sqrt{3}/4\}.$$
 
\begin{center}  \includegraphics[scale=0.4]{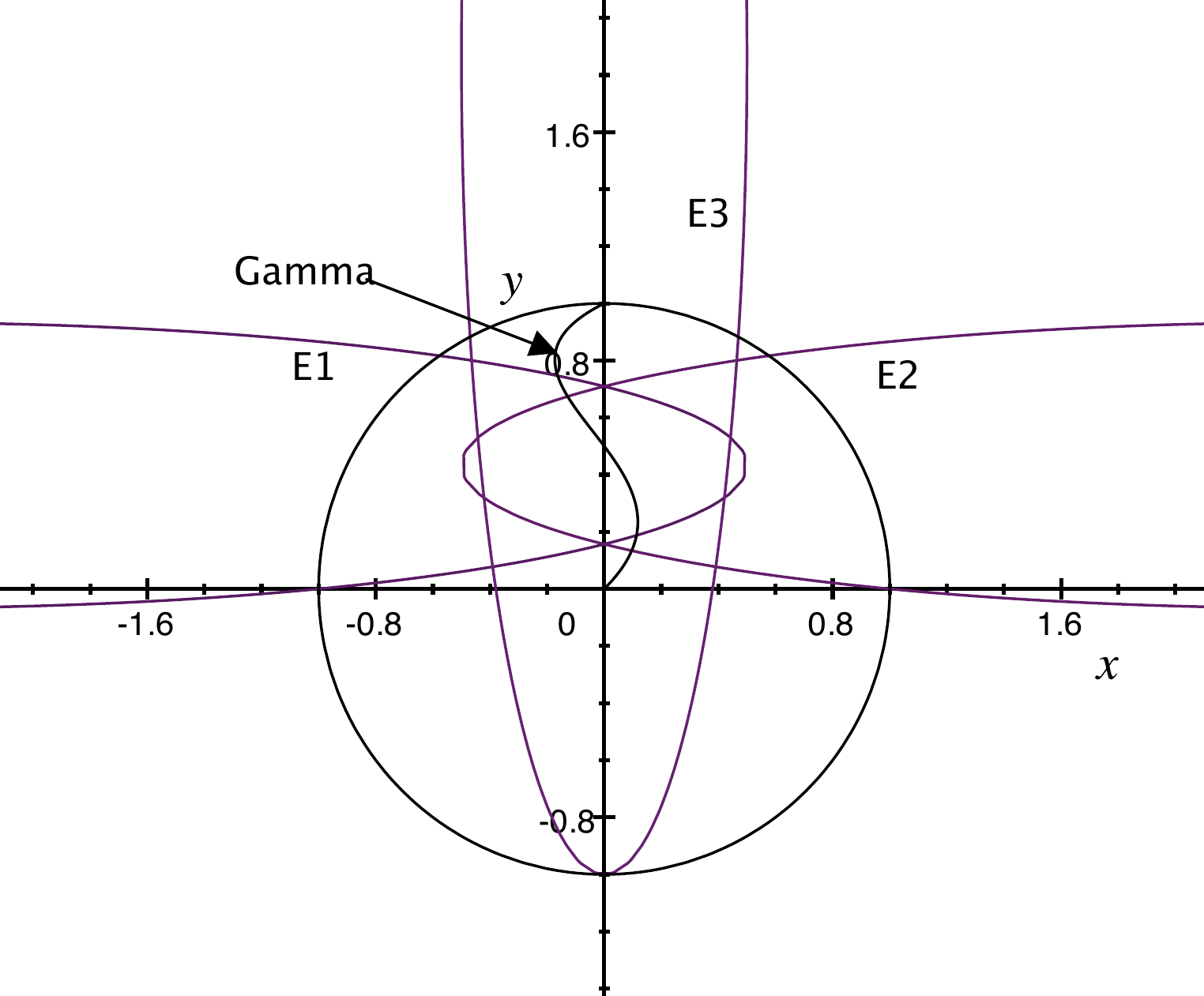}\end{center}
  
Let  $v_i$ be given by Lemma \ref{lemme1}, such that $0\leq v_i\leq 1$, $v_i = 0$ on the elliptic boundary of $E_i$  and  $v_i$  satisfies $F[v_i]+h(x)\cdot \nabla v_i|\nabla v_i|^\alpha + (V(x))v_i^{\alpha +1} > 0$. 
    
     There exists $\Gamma$ some  simple and regular curve which is included in $G$ and  links $(0,0)$ to $(0,1)$. 
    
    Let $E = E_1 \cup E_2$. We denote by $\partial E^+$ and $\partial E^-$ the superior and inferior boundary of $E$. Necessarily $\Gamma$  cuts $\partial E^+$ and $\partial E^-$. Let $\varphi$ be a parametrisation of $\Gamma$ with $\varphi$ in ${\cal C}^2$,  $\varphi(0) = (0,0)$ and $\varphi (1) = (0,1)$. Let 
    $t^-= \sup \{ t, \varphi (t)\in \partial E^-\}$ and 
    $t^+ = \inf \{ t, \varphi (t)\in \partial E^+\}$, and let 
    $p ^-= \varphi (t^-)$, $p^+ = \varphi (t^+)$. The portion of curve $(p^-, p^+)$ in $\Gamma$ is such that 
      for all $t \in ]t^-, t^+[$,  $\varphi(t) $ is in the interior of $E$. 
       Using the orientation of the portion  of curve between $p^-$ and $p^+$ one gets that  this portion of curve separates $E$ in two parts,  the left $E_l$ and the right $E_r$. 
       
        Let $z\in E_1\cap E_2$; if $z\in E_l$, we choose $D = E_2\cap E_l$, otherwise  $z\in E_r$ and $D = E_1 \cap E_r$. In  the first case $D$ has a boundary made of  parts of $\partial E_2$ and  the arc $\widehat{p^-,p^+}\cap E_2$. In the second one the boundary of $D$ has a boundary made of parts of $E_1$ and  $\widehat{ p^-p^+}  \cap E_1$.

 For example in the second case
 $$u-\chi_1v_1 >  \chi_1 (1-v_1)>0\  {\rm on } \ \widehat{p^-,p^+}\cap E_1$$
 $$ u-\chi_1 v_1 = u >0\ {\rm on} \ \partial  E_1$$
 and analogous  inequalities in the first case.

 Using the comparison principle in Theorem \ref{thcomp2}, we have obtained that 
 $$u(P) \geq \chi_1 \min \{v_1(
 P), v_2(P)\}\quad\mbox{ for all}\quad P\in E_1\cap E_2.$$ 
 
 Now we will use this to prove a similar inequality in $E_3$.

 One has   :  
  $$u\geq  \inf_{\{P\in \partial E_3, x_2= {\sqrt{3}/4}\}}  \min( v_1(P), v_2(P)) \chi_1 v_3.$$
  Indeed, this inequality holds on  $\partial E_3$, because on the elliptic part of $E_3$, $v_3=0$ and the straight part is included  in  $E_1\cap E_2$, where the inequality holds. Now by the comparison principle (Theorem \ref{thcomp2}) the inequality holds  true in $E_3$.
  
  We apply this in the ball $B(0,1/3)$ which is strictly included in the interior of $E_3$; defining 
  $$m_3 = \inf_{ B(0,1/3)} v_3,$$
 we have obtained that 
  \begin{eqnarray*}
  u&\geq & \chi_1\inf_{\{P\in \partial E_3, x_2= {\sqrt{3}/4}\}}(\min( v_1(P), v_2(P)) )m_3\\
  &\geq& {u(0)}\kappa\inf_{\{P\in \partial E_3, x_2= {\sqrt{3}/4}\}}(\min( v_1(P), v_2(P)) )m_3.
  \end{eqnarray*}

\bigskip

\noindent{\em Proof of Claim 2}. Fix any point $\overline P$ in $B_{\frac{1}{4}}(0)$. Then
$$ 
 B_{\frac{3}{4}}(\overline P)\subset B_1(0), \ \mbox{and}  \ 0\in  B_{\frac{1}{4}}(\overline P). $$
 Hence by Claim 1 we have that
 $$u(\overline P)\leq Ku(0)$$
 but always by Claim 1,  $u(0)\leq Ku(\overline P)$. 
 This ends the proof of Claim 2, by choosing $K_1=\frac{1}{K}$ and $K_2=K$.
 
 \bigskip
 \noindent{\em Proof of Claim 3}. Now $\Omega=B(0,R)$ and $u$ is a positive solution of (\ref{eqhar}).
Using the homogeneity of $F$, let
$v(p):=u(Rp)$ satisfies
$$F(x, \nabla v, D^2 v)+ R h(R x)\cdot \nabla v|\nabla v|^\alpha +  R^{\alpha+2}( V( Rx)) v^{\alpha+1}\leq0, \mbox{ in}\ B_{1}(0).$$
Hence we are in the  conditions of the previous case with $h$ replaced by $Rh(Rx)$ and $V(x)$ replaced by $R^{\alpha+2}V(Rx)$.  We have obtained that $v$ satisfies, for any $P\in B_{\frac{R}{3}}(0)$:
$$v(0)\leq Kv(P)\ \mbox{i.e.} \ u(0)\leq Ku(Q)\ \mbox{for} \ Q\in  B_{\frac{R}{3}}(0).$$
Observe that $K$ depends on $\gamma$ (see (\ref{gam})), but when $R\leq R_o$ it can be chosen independently on $R$. (Moreover let us note that in Liouville's result we shall consider arbitrary large values  of $R$, but since we shall assume that $h = V= 0$, $\gamma $ is independent of $R$ and this will allow to get the result. )

 \bigskip
  \noindent{\em Proof of Claim 4}. This is standard potential theory procedure. Let $K$ be a compact connected subset of $\Omega$. And let $R=\inf\{r,  d(P,\partial\Omega)\leq r, \mbox{ for any } P\in K \}.$ 
  Suppose that $P$ and $Q$ are any two points of $K$. Then there exists a continuous curve $\Gamma\subset K$ joining $P$ and $Q$. We can find a finite number of points $P=P_1,P_2,\dots,P_k=Q$ such that 
  $$P_i\in\Gamma,\ |P_iP_{i+1}|\leq\frac{R}{4},\ B_R(P_i)\subset \Omega.$$
 Hence applying the previous results, observing that  $$\overline\lambda(\Omega)<\overline\lambda (B_R(P_i))$$
we get
$$u(P)\leq K_2u(P_2)\leq K_2^2u(P_3)\leq K_2^ku(Q).$$
This ends the proof of  Theorem \ref{prop1}.

\bigskip
 
Proof of Theorem \ref{prop2}

We recall that $V = 0$ and $f\neq 0$  and we shall give  shortly the changes in the case $V \leq 0$ to prove the result in remark \ref{remV}. 
The proof proceeds with the same steps as in the case $f\equiv 0$, the difference being that instead of comparing $u$ with the functions $v_i$ defined in Lemma \ref{lemme1} we will need to compare it with $Cv_i+w$ where $w$ is a subsolution of (\ref{eqhar2}), and since the operator is fully-nonlinear we need to prove that  $Cv_i+w$ is a subsolution.

We begin with the case $\alpha \leq 0$. 
We define  $q = {\alpha + 2\over \alpha +1}$.

Let 
$\rho_o$  and $C_1$ such that 
$$\rho_o = \inf (1,{1\over 8 |h|_\infty  aq^2 }),\quad
\mbox{and}\quad C_1 =({2|f|_\infty  \over aq^{2+\alpha}})^{1\over 1+\alpha}.$$
Later we shall also impose to $C_1$ to be greater than some given constants. 
As in  the case $f = 0$, we need to prove that there exists some point on the boundary of $B(0,\rho_o)$ where  $u > {u(0)\over 2}-C_1\rho^q$. This is a consequence of the fact that with this choice of $C_1$ and $\rho_o$,  ${u(0)\over 2}- C_1 \rho^q$ is a supersolution  in the set $\rho < \rho_o$. 
We can assume that this point is $(0,\rho_o)$, and we define the curve $\Gamma$ as in the case where $f\equiv 0$ 

We introduce the functions $v_i$ for $i = 1,2,3 $  related to the ellipsis $E_i$ defined as 

          $$E_1 = \{ (x_1,x_2),\ \sigma_1^2:={(x_1+{5\rho_o\over 2})^2\over 9\rho_o^2 } + 4\frac{(x_2-{\rho_o\sqrt{3}\over 4})^2}{\rho_o^2} \leq 1, \ x_1\geq -\rho_o\}$$
  $$E_2 = \{ (x_1,x_2),\ \sigma_2^2:= {(x_1-\rho_o{5\over 2})^2\over 9\rho_o^2 } + 4\frac{(x_2-{\rho_o\sqrt{3}\over 4})^2}{\rho_o^2} \leq 1, \ x_1\leq \rho_o \}$$
 
 $$E_3 = \{ (x_1,x_2), \ \sigma_3^2:=\frac{4x_1^2}{\rho_o^2}+ {\left(x_2- \rho_o(1+{\sqrt{3}\over 2})\right)^2\over \rho_o^2(2+ {\sqrt{3}\over 2})^2} \leq 1, x_2\leq \frac{\sqrt{3}\rho_o}{4}\}.$$
 
 Recall that 
 $$v_i = {e^{-\gamma_i \sigma_i^2}-e^{-\gamma_i}\over e^{-\gamma_i\over 4}-e^{-\gamma_i}}$$ where  e.g. 
 $$ \gamma_1 =\sup \left({4(A+a)\over a} (1+  36), {4|h|_\infty M_19\over m_1a} \right)$$
 for
 $$m _1= 3^{-\alpha}, \quad\mbox{and}\quad M_1 = 2^{1+\alpha} ({1\over 9} + 4)^{1+\alpha\over 2}.$$

 For the following  we shall replace the constant $\gamma_i $ by $ \gamma\equiv \sup \gamma_i$  which is also convenient  to our  goal. 

 We need to observe that 
 $$|\nabla v_i |\leq {4\gamma\over \rho_o} \tilde v $$
   where 
  $\tilde v= {e^{-\gamma\sigma^2}\over e^{-\gamma\over 4}-e^{-\gamma}}$.
  Note that $v_o= {e^{-\gamma\over 4}\over e^{-\gamma\over 4}-e^{-\gamma}}\geq \tilde v\geq e^{-3\gamma\over 4} v_o $. 
  With all these choices of constants, the computation in Lemma \ref{lemme1} gives for $i = 1,2,3 $
  \begin{eqnarray*}
  |\nabla v_i|^\alpha ({\cal M}^-(D^2v)-h(x)\cdot \nabla v_i)
   &\geq& {2^{2\alpha-1} \gamma^{2+\alpha} a \tilde v^{1+\alpha} \over 9 \rho_o^{2+\alpha}}\\
    &\geq & {2^{2\alpha-1} \gamma^{2+\alpha} a  (e^{-3\gamma\over 4} v_o)^{1+\alpha} \over 9 \rho_o^{2+\alpha}}\\
       & :=& {c_2\over \rho_o^{\alpha+2}}
   \end{eqnarray*}

   We now consider two cases  : 
   
  - Either $\displaystyle{{u(0)^{1+\alpha} c_2\over 2^{1+\alpha}\rho_o^{\alpha+2}}}> |f|_\infty$
    and then  for $i=1,2,3$ ${u(0) v_i\over 2}$ is a subsolution of  the equation, while $u+2C_1\rho_o^q$ is a supersolution of the same equation. 
     Moreover in $E_1\cap E_2$, using the fact that the boundary is made of arcs of $\Gamma$ or of parts of the boundary of $E_i$ one gets that 
     
     $u+ 2C_1\rho_o^q \geq {u(0)\over 2} \inf (v_1, v_2)$ in $E_1\cap E_2$. 
     
     \noindent And  the final step is as in the case  where $f=0$, i.e.  proving that 
      $u+ 2C_1\rho_o^q \geq {u(0)\over 2} \inf_{x_2 = {\sqrt{3}\over 4} , x\in \partial E_3} (v_1, v_2) v_3$

 - Or $\displaystyle{{u(0)^{1+\alpha} c_2\over 2^{1+\alpha}\rho_o^{\alpha+2}}}\leq  |f|_\infty$ 
  
In that case let  $\tilde \rho = (x_1^2+ (x_2+ 3\rho_o)^2)^{1\over 2}$.     Let us note that $\rho\leq \rho_o \leq {\tilde \rho\over 2} \leq \tilde \rho \leq 4\rho_o$.  In particular on the curve $\Gamma$, $u \geq {u(0)\over 2}- C_1\rho^q\geq {u(0)\over 2}-C_1\tilde \rho^q$. 
 
 We choose $C_1 = \sup 
 \left({e^{3\gamma\over 4}9^{1\over 1+\alpha}
  \over (a\gamma)^{1\over 1+\alpha}}  ,  \left({2^{1-\alpha} \over aq^{2+\alpha}}\right)^{1\over 1+\alpha} \right) |f|_\infty^{1\over 1+\alpha}$ 
 $$w =   C_1\tilde  
 \rho^q,$$

 We shall prove that  for $i=1,2$
 \begin{equation}\label{eq1}
 (|{u(0)\over 2}\nabla v_i|+ |\nabla w|)^\alpha \left( {\cal M}^- (D^2 w)-|h|_\infty |\nabla w|)\right) \geq |f|_\infty
 \end{equation}
 in $E_1\cap E_2$  
  and 
   \begin{equation}\label{eq2}
  (|{u(0 )\over 2}\nabla v_3|+ |\nabla w|)^\alpha  \left({\cal M}^- (D^2 w)-|h|_\infty |\nabla w|\right) \geq |f|_\infty
  \end{equation}
  in $E_3$. 
  
  For that aim  we observe that $|{u(0)\over 2}\nabla v_i|\leq |\nabla w|$ by the choice of $C_1$. 
  For simplicity we shall do the computation only for   $v_1$. Observe first that
  \begin{eqnarray*}
  |{u(0)\over 2}\nabla v_1|&\leq& {2 u(0)\gamma \tilde v \over \rho_o} \\
  &\leq& 2 ({u(0) v_o\gamma \over \rho_o}) \\
  &\leq &e^{3\gamma\over 4} \left( {|f|_\infty \rho_o 9 \over a\gamma 2^{\alpha-2}} \right)^{1\over 1+\alpha}\\
&  \leq& C_1 q (2\rho_o)^{q-1} \leq qC_1 \tilde \rho^{q-1}.
\end{eqnarray*}
 From this and similar calculations, we derive that,  for $i = 1,2,3 $,
     \begin{eqnarray}
    (|\nabla \left({u(0)\over 2}v_i\right)|+ |\nabla w|)^\alpha \left( {\cal M}^- (D^2 w)-|h|_\infty |\nabla w|\right) &\geq&
     (2qC_1 \tilde \rho^{q-1})^\alpha {aq^2\over 2} C_1 \tilde \rho^{q-2} \nonumber \\
   & \geq & 2^\alpha C_1^{1+\alpha} q^\alpha {aq^2\over 2} \nonumber \\
   &\geq& |f|_\infty .\label{eq4}
   \end{eqnarray}
      
Moreover  from the choice of $\gamma_i$, one has  
$${\cal M^-} (D^2 v_i)-|h|_\infty |\nabla v_i|  \geq 0,$$  and, using the simple inequality 
$|X+ Y|^\alpha\geq |X|^\alpha + |Y|^\alpha$,  this implies that  ${u(0)v_i\over 2} + w$ is a subsolution  for $i=1,2$ of the equation 
$$F[{u(0) v_i\over 2} + w] -h\cdot \nabla ({u(0) v_i\over 2} + w)|\nabla ({u(0)  v_i\over 2} + w)|^\alpha \geq f\quad\mbox{in} \quad E_1\cap E_2.$$
 We have obtained that  
 \begin{equation}\label{040209}u+ 2(4\rho_o)^q C_1 \geq {u(0) \inf_{i=1,2}v_i\over 2} + w\quad\mbox{in} \quad E_1\cap E_2
 \end{equation}
  if it is true on the boundary of $E_1\cap E_2$.
On the elliptic boundary  of $E_1$  since $u>0$ it is immediate, while    
on the part of the boundary made of   portions of $\Gamma$

  $$u\geq{ u(0)\over 2} -C_1\rho^q \geq {u(0)\over 2} -C_1\tilde \rho^q \geq{ u(0)\over 2}-(4\rho_o)^q C_1$$
   and then
    (\ref{040209}) holds true.

Finally we remark that $\chi= \inf_{\{x\in \partial E_3,\  x_2 = {\sqrt{3}\over 4} , i=1,2\}}\inf (v_i)\leq 1$  and then from the 
 equation (\ref{eq4})
$$F[{u(0)\over 2} \chi v_3 + w] -h\cdot \nabla ({u(0) \chi v_3 \over 2}+ w)|\nabla ({u(0) \chi v_3\over 2} + w)|^\alpha \geq f. $$
 ${u(0)\over 2} \chi v_3 + w$ is  then a sub-solution, which satisfies  on the boundary of $E_3$ the inequality 
   $u+ 2(4\rho_o)^q C_1\geq {u(0)\chi \over 2}v_3+ w$,  since this is  true   on the straight part  of $E_3$ which is included in   $\{x_2 = {\sqrt{3}\over 4},x\in \partial E_3\}$ and it is true on the elliptic part of $E_3$ because $v_3=0$  on this part and $u>0$.
    We have obtained that there exist some constant $K$ and $K^\prime$ such that 
 $$u \geq    K u(0) -K^\prime |f|_\infty^{1\over 1+\alpha}.$$
    This ends the proof of the case $\alpha \leq 0$ and $V = 0$. 
    \bigskip

           \bigskip
        
        We now consider the case $\alpha >0$ and $V=0$. 
The first part of the proof proceeds as for the case $f=0$ : For a fixed $\delta\in (0,1)$ that will be introduced later, we  define  : 
         
      $\rho_o = \inf ({2^{-{|\alpha-2|\over 2}-\alpha-1} \delta^\alpha aq \over  |h|_\infty},1)$,   $C_1= \left( {2^{{|\alpha-2|\over 2}+ 2} |f|_\infty \over a q^{2+\alpha} \delta^\alpha  }\right)^{1\over 1+\alpha}$, 
        and  the function $w_1 = C_1 \rho^q$,   it is clear that  $ w_1 $ satisfies
          $$F[-w_1] + h\cdot \nabla (-w_1)|\nabla w_1|^\alpha \leq -|f|_\infty, $$
 then so does ${u(0)\over 2}- w_1.$        
 
  Let $G_1 = \{x\in B(0,\rho_o), u(x) > {u(0)\over 2}-C_1 \rho^q\}$.
   $G_1$ is an open set which contains $0$.  Let $G$ be the connected component of $G_1$ which contains $0$.  By the comparison principle the boundary of  $G$ contains at least one point of $\partial B(0, \rho_o)$.  One can assume that this point is $(0, \rho_o)$.

   Let $\tilde\rho^2=x_1^2+(x_2+3\rho_o)^2$ and $w=C_1\tilde \rho^q$. Let $\Gamma$ be a regular curve which links $0$ to $(0, \rho_o)$ and is included in $G$, then  since $\tilde \rho > \rho$, one always has  $u > {u(0)\over 2}-C_1\tilde \rho^q$ on $\Gamma$.

         We now proceed to the second step. 
 From Lemma \ref{lemme1}, with $\sigma:= \sigma_i$ associated to the ellipsis $E_i$ we know that the function $v_i= {e^{-\gamma_i\sigma_i^2}- e^{-\gamma_i}\over e^{-\gamma_i\over 4}- e^{\-\gamma_i}}$ is a subsolution  of:

$$\delta^\alpha 2^{-|\alpha -2|\over 2} |\nabla v_i|^\alpha {\cal M}_{a,A}^- (D^2 v_i) -2^\alpha  |h|_\infty |\nabla v_i|^{\alpha+1}\geq 0$$
in $E_i$, with an appropriate choice of $\gamma_i$  e.g. 
$$\gamma_1 = \sup \left({4(A+a)\over a}(1+ 9/4) , {2^{{|\alpha-2|\over 2} + \alpha+2} 9 |h|_\infty M_1\over \delta^\alpha m_1 a}\right)$$
 for some obvious definitions of $m_1$ and $M_1$. 

For $i=1$ and $i=2$, we need to show that in $E_1\cap E_2$, ${u(0)\over 2} v_i+w$ is a subsolution of (\ref{eqhar}).
To do so we need to evaluate $\grad v_i\cdot\grad w$; this is done in Lemma \ref{lem2}  below. Applying it,  there exists some $1>\delta>0$ such that 
 for $i = 1,2$ one has in $E_1\cap E_2$, 
\begin{eqnarray*}
         |\nabla \left({ u(0)\over 2} v_i\right)+ \nabla w|^2& \geq & |\nabla \left({ u(0)\over 2} v_i\right)|^2+ |\nabla w|^2  + 2(-1+ \delta^2)|\nabla\left( { u(0)\over 2}v_i\right) ||\nabla w|\\
         &\geq& \delta^2 |\nabla\left( { u(0)\over 2} v_i\right)|^2 + \delta^2 |\nabla w|^2
 \end{eqnarray*}
  and  in $E_3$ denoting as $\chi$ the constant $\chi = \inf_{\{P\in \partial E_3, x_2= {\sqrt{3}/4}\}}  \min( v_1(P), v_2(P))$
       $$|\nabla\left({ u(0)\over 2}\chi  v_3 \right)+ \nabla w|^2
         \geq \delta^2 |\nabla \left({ u(0)\chi\over 2}   v_3\right)|^2 + \delta^2 |\nabla w|^2.
         $$
         
            Let us note that with the choice of $C_1$ made in the first step, one has 
    
          $2^{-|\alpha-2|\over 2}\delta^\alpha  |\nabla w|^\alpha {\cal M}^- (D^2 w)-2^\alpha |h|_\infty |\nabla w|^{1+\alpha }  \geq aq^{2+\alpha}  2^{{-|\alpha-2|\over 2}-2}\delta ^{\alpha} C_1^{1+\alpha}\geq    |f|_\infty$. 

 This implies that  for $i=1,2$, in $E_i$
         \begin{eqnarray*}
  &&        |{u(0)\over 2}\nabla v_i+ \nabla w|^{\alpha}{\cal M}^- ({u(0)\over 2} D^2v_i+ D^2 w)
          -
 |{u(0)\over 2} \nabla v_i + \nabla w|^{\alpha+1}|h|_\infty \\
      &\geq &    2^{-|\alpha-2|\over 2} \delta^\alpha( |{u(0)\over 2} \nabla v_i|^\alpha  + |\nabla w|^{\alpha}) {\cal M} ^-({u(0)\over 2} D^2v_i+ D^2 w)\\
      &&-
          2^\alpha 
          |{u(0)\over 2} (\nabla v_i|^{\alpha+1}  + |\nabla w|^{1+\alpha}) \\
&\geq&         2^{-{|\alpha-2|\over 2}} \delta ^\alpha |{u(0)\over 2} \nabla v_i|^\alpha {\cal M}^- ( {u(0)\over 2} D^2v_i )\\
 &&-        |h|_\infty 2^{\alpha } |{u(0)\over 2} \nabla v_i|^{1+\alpha} +         2^{-|\alpha-2|\over 2} \delta ^\alpha |\nabla w|^\alpha {\cal M}^- (D^2 w) -|h|_\infty 2^{\alpha } |\nabla w|^{1+\alpha}\geq |f|_\infty
 \end{eqnarray*}          
         and also 
         \begin{eqnarray*}
  &&        |{u(0)\over 2}\chi\nabla v_3 + \nabla w|^{\alpha} {\cal M}^- ({u(0)\over 2}\chi  D^2v_3+ D^2 w)
          -
          |{u(0)\over 2}\chi \nabla v_3 + \nabla w|^{\alpha+1}|h|_\infty \\
                &\geq & 2^{-|\alpha-2|\over 2} \delta ^\alpha (|{u(0)\over 2}\chi \nabla v_3|^\alpha {\cal M}^- ( {u(0)\over 2} \chi D^2v_3 )\\
          &&-|h|_\infty 2^{\alpha } |{u(0)\over 2}\chi \nabla v_3|^{1+\alpha} + 2^{-|\alpha-2|\over 2} \delta ^\alpha |\nabla w|^\alpha {\cal M}^- (D^2 w) -|h|_\infty 2^{\alpha } |\nabla w|^{1+\alpha}\geq |f|_\infty.
\end{eqnarray*}

    Let $\tilde \rho_o= 4 \rho_o$ in order that in $\rho< \rho_o$, $\tilde \rho < \tilde \rho_o$. 
     We check that 
     $u + 2C_1\tilde  \rho_o^q \geq {u(0)\over 2} \inf_{i=1,2}(v_i) + w$ in $E_1\cap E_2$. Indeed,  $u+2C_1\tilde  \rho_o^q $ is a supersolution and      ${u(0)\over 2} v_i + w$ is a sub-solution of the same equation in $E_i$, $i = 1,2$. 
      Moreover  if $x\in D$ which is made of some part of $\partial E_1$ and some part of $\Gamma$ one gets that  $u+2C_1\tilde  \rho_o^q >  w$ on  the boundary of $\partial E_1$ since   $u$ is positive. On $\Gamma$ it is true since 
      $u\geq{u(0)\over 2} -C_1\tilde \rho^q $.
      
     We now proceed to the last part of the proof : 
      
      We have on the straight part of $E_3$
      $u+ 2 C_1  \tilde\rho_o^q\geq {u(0) \chi v_3\over 2} + C_1 \tilde \rho^q$. 
       Indeed one has on that part 
       $u \geq     \inf_{\{P\in \partial E_3, x_2= {\sqrt{3}/4}\}}  \min( v_1(P), v_2(P)){u(0) v_3\over 2} - C_1  \tilde \rho^q$ . 
       
        On the elliptic part of $E_3$ the result is true since $ v_3 = 0$ . Since 
        $ {u(0) \chi v_3\over 2} + C_1 \tilde \rho^q$ is a sub-solution and $u+ 2 C_1  \tilde\rho_o^q$ is a supersolution we have obtained the result, as in the case $f = 0$ that on $B(0, \rho_o/3)$ there exists some constant $K$ and $K^\prime$  which do not depend on $\rho_o<1$, such that 
        $$ u\geq Ku(0)-K^\prime|f|_\infty^{1\over 1+\alpha}$$
        The rest of the proof is the same.

  \begin{lemme}\label{lem1}
  Let $\rho <1$. 
   Let $w_1 =  C_1 \tilde \rho^q$. Then for $ C_1 = \left({|f|_\infty 2^{{|\alpha-2|\over 2}+ 1} \over \delta^\alpha a q^{2+\alpha}}\right)^{1\over 1+\alpha}$, 
   $$\delta 2^{-|\alpha-2|\over 2} |\nabla w|^\alpha {\cal M}^- (D^2 w) -|h|_\infty 2^{\alpha} |\nabla w|^{1+\alpha} \geq |f|_\infty$$ in the set $E_3$. 
     \end{lemme}
   
   Proof:

   One has  for $Cx = (x_1, x_2+3\rho_o)$, 
   $\nabla  w_1 =  C_1q \tilde \rho^{q-2} Cx$
   and $\nabla \nabla w =  C_1q \tilde\rho^{q-4} ((q-2) Cx\otimes  Cx+\tilde  \rho^2 I)$. 
   The matrix $ (q-2) Cx\otimes  Cx+ \tilde \rho^2 I$ has eigenvalue $(q-1) \tilde \rho^2$ and $\tilde \rho^2$, as a consequence 
   $$\delta^\alpha 2^{-|\alpha-2|\over 2}|\nabla w|^\alpha {\cal M}^- (D^2 w) -|h|_\infty |\nabla w|^{\alpha+1} 2^\alpha \geq \delta^\alpha 2^{-|\alpha-2|\over 2} aq^{2+\alpha}\tilde \rho^{(q-1)\alpha + q-2} -| h|_\infty \tilde \rho 2^{\alpha} \geq 
 |f|_\infty.$$

\begin{lemme}\label{lem2}
There exists $\delta\in [0,1[$ such that in $E_1\cap E_2$ for $i=1$ and $i=2$
$$\langle \nabla v_i,\nabla w\rangle    \geq (-1+\delta^2) |\nabla v_i| |\nabla w|$$
 and in $E_3$
 $$\langle \nabla v_3,\nabla w\rangle \geq (-1+\delta^2)  |\nabla v_3|\ |\nabla w|.$$
 
\end{lemme}
{\em Proof }
For homogeneity reasons, we can assume that $\rho_o = 1$. Then $\nabla v_i =\gamma_i  B_i x \tilde v $, with 
  $B_1 x :=-({x_1+{5\over 2}\over 9} , 4(x_2-{\sqrt{3}\over 4}))$, 
          $B_2 x := - ({x_1-{5\over 2}\over 9} , 4(x_2-{\sqrt{3}\over 4}))$,  and $B_3 x=-(4x_1, {x_2-1-{\sqrt{3}\over 2}\over (2+ {\sqrt{3}\over 2})^2})$. While $\grad w=C_1\tilde\rho^{q-2}(Cx)$ with 
$Cx = (x_1, x_2+3)$.  

           It is an elementary but tedious calculation to see that for $x\in E_1\cap E_2$ the vectors $B_1x$, $B_2x$  lie in the circular sector $S$ defined by 
$\frac{6\sqrt{11}}{5}|x_1|\geq x_2\geq\frac{-6\sqrt{11}}{5}|x_1|$, while $Cx$ lies in a sector $S_o$
defined by $\frac{\sqrt{3}+12}{2}|x_1|\leq x_2$. Hence if $\theta_1$ is the angle between the sectors then the first equality is satisfied with $-1+\delta=\cos\theta_1$.
Similarly for the second case.

\begin{center} \includegraphics[scale=0.3]{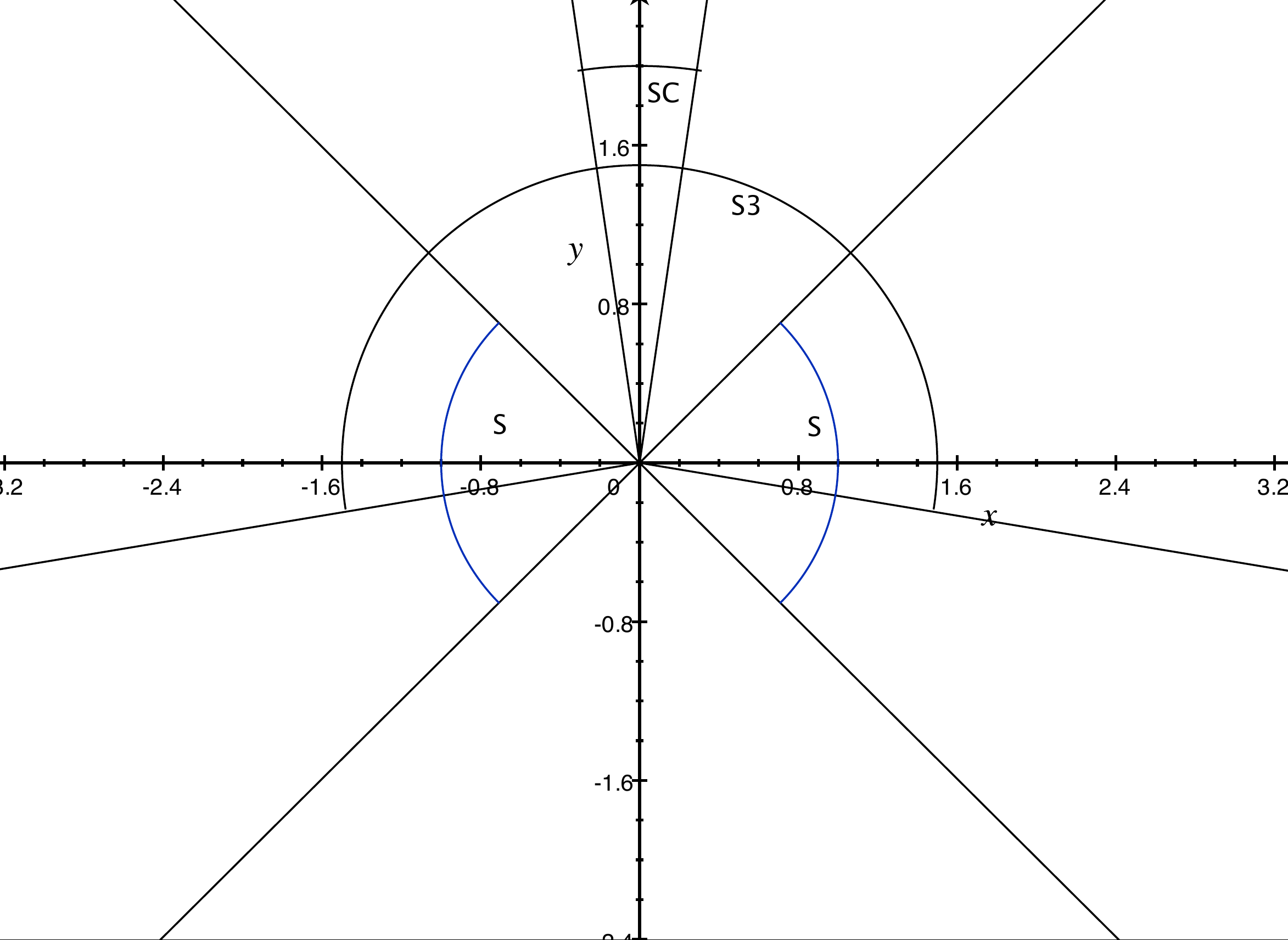}\end{center}
 
 The  circles of smaller radius indicate  the sectors spanned by $B_ix$ and the circle of larger radius indicates the sector spanned by $Cx$, as can be seen the angle between $B_ix$ and $Cx$ is never $\pi$.

  \bigskip
\noindent {\bf The case $V\leq 0$, $f\neq 0$}

As in the previous proof we begin with the case $\alpha <0$. 

  We now consider the case $V\leq 0$. We choose $\rho_o = \inf \left({aq\over 4 |h|_\infty }, 
     \left({aq^{2+\alpha} \over |V|_\infty 8}\right)^{1\over \alpha+2}\right)$ 
     and 
     $C_1=  \left({2^{1-\alpha} \over aq^{2+\alpha}}\right)^{1\over 1+\alpha}$ and we choose the constant $\gamma_i$  as in Lemma \ref{lemme1} with $V\leq 0$, in order that 
     $|\nabla v_i|^\alpha {\cal M}^-(D^2 v_i)-|h|_\infty |\nabla v_i|^{1+\alpha}-|V|_\infty v_i^{1+\alpha}) \geq 0$. 
     
      Let us note that since $V\leq 0$, $u$ is  also a subsolution of 
      $$|\nabla u|^\alpha\left ({\cal M}^-(D^2 u)-h(x)\cdot \nabla u\right) \geq f$$
       and then the first step is still valid.  We obtain that there exists some point on the boundary  of $\partial B(0, \rho_o)$ such that 
       $u \geq {u(0)\over 2} - C_1 \rho_o^q$.  We can assume that this point is $(0, \rho_o)$.

      We now consider as previously two cases:
      
 -Either $\displaystyle{{u(0)^{1+\alpha} c_2\over 2^{1+\alpha}\rho_o^{\alpha+2}}}> |f|_\infty$
    and then  for $i=1,2$ ${u(0) v_i\over 2}$ is a subsolution of  the equation in $E_1\cap E_2$ and   so is ${u(0) \chi v_3\over 2}$ in $E_3$ , while $u+2C_1\rho_o^q$ is a supersolution of the same equation. 
     Moreover in $E_1\cap E_2$, using the fact that the boundary is made of arcs of $\Gamma$ or of parts of the boundary of $E_i$ one gets that 
     $u+ 2C_1(4\rho_o)^q \geq {u(0)\over 2} \inf (v_1, v_2)$
in $E_1\cap E_2$. 
And now we do the final step as in the case  where $f=0$, i.e. we prove that 
      $u+ 2C_1(4\rho_o)^q \geq {u(0)\over 2} \inf_{x_2 = {\sqrt{3}\over 4} , x\in \partial E_3} (v_1, v_2) v_3$.

 - Or $\displaystyle{{u(0)^{1+\alpha} c_2\over 2^{1+\alpha}\rho_o^{\alpha+2}}}\leq  |f|_\infty$. 
 
 In that case we  define $C_1=  \sup \left(e^{3\gamma\over 4} {\left(9
  \over a\gamma \right)^{1\over 1+\alpha}} ,  \left({2^{1-\alpha} \over aq^{2+\alpha}}\right)^{1\over 1+\alpha},  2^{1\over 1+\alpha} c_2^{-1\over 1+\alpha} \right) |f|_\infty^{1\over 1+\alpha}$   and $\tilde \rho$ as in the case $V=0$ and we observe that 
  $|\nabla v_i |{u(0)\over 2}\leq |\nabla w|$. Let us note that, with the choice  of $C_1$ above,   $${v_i u(0)\over 2} \leq w$$
  and 
  $$2^\alpha |\nabla w|^\alpha\left( {\cal M}^-(D^2 w)-|h|_\infty |\nabla w|\right)-2^{\alpha +1}|V|_\infty  w^{1+\alpha} \geq |f|_\infty.$$
  We now write  for $i = 1,2$ 
  \begin{eqnarray*}
&& ( |\nabla {u(0) v_i\over 2}|+ |\nabla w|)^\alpha\left({\cal M}^- (D^2 ({u(0)v_i\over 2} + w) -h(x)\cdot (\nabla {u(0) v_i\over 2} + \nabla w)\right)\\
 &-& |V|_\infty  ({u(0)v_i\over 2} + w) ^{1+\alpha} \\
 &\geq & 
   ( |\nabla {u(0) v_i\over 2}|+ |\nabla w|)^\alpha\left({\cal M}^- (D^2 ({u(0)v_i\over 2})) -h(x)\cdot (\nabla {u(0) v_i\over 2}) \right)\\
   &+ &  (2 |\nabla w|)^\alpha\left({\cal M}^- (D^2 ( w) -h(x)\cdot \nabla w\right)-2^{\alpha+1}|V|_\infty  w^{1+\alpha}\\&\geq & 0+ |f|_\infty
   \end{eqnarray*}
   
   and for $i = 3$ and $\chi = \inf_{\{x_2 = {\sqrt{3}\over 4}, x\in \partial E_3\}} \inf(v_1, v_2)$
     \begin{eqnarray*}
&& ( |\nabla {u(0)\chi v_3\over 2}|+ |\nabla w|)^\alpha\left({\cal M}^- \left((D^2 ({u(0)\chi v_3\over 2} + w) \right)-h(x)\cdot (\nabla {u(0) \chi v_3\over 2} + \nabla w)\right)\\
 &-& |V|_\infty  ({u(0)\chi v_3\over 2} + w) ^{1+\alpha} \\
 &\geq & 
   ( |\nabla {u(0)\chi  v_3\over 2}|+ |\nabla w|)^\alpha\left({\cal M}^- (D^2 ({u(0)\chi v_3\over 2})) -h(x)\cdot (\nabla {u(0) \chi v_3\over 2} )\right)\\
   &+ &  (2 |\nabla w|)^\alpha\left({\cal M}^- (D^2 ( w)) -h(x)\cdot \nabla w\right)-2^{\alpha+1}|V|_\infty  w^{1+\alpha}\\&\geq & |f|_\infty.
   \end{eqnarray*}

   The rest of the proof is analogous to the one done in the previous cases, observing that, since $V \leq 0$, $u+ 2C_1(4\rho_o)^q$ is also a supersolution of the equation.

\bigskip
 
  We now  treat the case $\alpha >0$.       The notations $B_i$,   $C$, $\delta$ are the same as in the  case $V=0$. 
      
     Since  $V\leq 0$, $u$ is also a subsolution of 
      $F[u] + h(x)\cdot \nabla u |\nabla u|^\alpha \geq f$ so the first step is the same, more precisely  if 
      we choose $\rho_o < \inf \left(\left({|h|_\infty \over 4|V|_\infty}\right)^{1+\alpha}, {\delta ^\alpha 2^{{-|\alpha-2|\over 2}-\alpha-2}aq^2 \over |h|_\infty}\right)$,  and  $C_1= \left( {2^{{|\alpha-2|\over 2}+ 3} |f|_\infty \over a q^{2+\alpha}  \delta^\alpha  }\right)^{1\over 1+\alpha}$, where $\delta$ is as in the proof of $V=0$, $\alpha >0$ and $f\neq 0$, 
$w_1 =- C_1 \rho^q$ is a supersolution of 
$F[w_1] +h(x)\cdot \nabla w_1 |\nabla w_1|^\alpha \leq -|f|_\infty$, then so is ${u(0)\over 2}+ w_1$. We obtain always by the some reasoning that there exists some point on the boundary $\rho = \rho_o$ on which  $u > {u(0)\over 2}- C_1 \rho^q$.

   For the second step we must prove that 
      one can chose $v_i$ such that  in $E_i$ 
           $$ \delta^\alpha 2^{-|\alpha-2|\over 2}{\cal M}^-(D^2 v_i)|\nabla v_i|^\alpha -2^{\alpha}|h|_\infty |\nabla v_i|^{\alpha+1}-2^{\alpha}|V|_\infty v_i^{1+\alpha} > 0. $$
           This can be done by choosing  $\gamma_i$ such that 
    $$\gamma_i = \sup \left({4(A+a)\over a} (1+  {b_i^2\over  c_i^2}), {2^{{|\alpha-2|\over 2} +\alpha+ 3}|h|_\infty M_ib_i^2\over \delta^\alpha  m_ia}, \left({2^{{|\alpha-2|\over 2}+ \alpha + 2} |V|_\infty b_i^2 \over am_i \delta^\alpha }\right)^{1\over 2+\alpha}\right)
$$
(with obvious definitions of $b_i,$ $c_i$, $M_i, m_i$, on the model of the proof of lemma \ref{lemme1}). 
    
     Let $\tilde \rho$ be defined as in the previous proof, then $w = C_1\tilde \rho^q$ is a solution of 
     
$2^{-|\alpha -2|\over 2} \delta^\alpha {\cal M}^- (D^2 w) -2^\alpha |h(x)|_\infty | \nabla w |^{\alpha+1}-2^\alpha |V|_\infty w^{1+\alpha } \geq |f|_\infty$,
 and then 
 ${u(0)\over 2} v_i + w$ is for $i=1,2$ a sub-solution of 
 $$F[{u(0)\over 2} v_i + w] + h(x)\cdot \nabla( {u(0)\over 2} v_i + w)|\nabla ({u(0)\over 2} v_i + w)|^\alpha
 + V(x)({u(0)\over 2} v_i + w)^{\alpha+1} \geq |f|_\infty$$
 in $E_1\cap E_2$ and  
         \begin{eqnarray*}
         F[{u(0)\over 2} \chi v_3 + w] &+& h(x)\cdot \nabla( {u(0)\over 2} \chi v_3 + w)|\nabla ({u(0)\over 2} \chi v_3 + w)|^\alpha+\\
 &+& V(x)\left( {u(0)\over 2} \chi v_3 + w\right)^{1+\alpha} 
 \geq |f|_\infty
 \end{eqnarray*}
in $E_3,$   with  $\chi = \inf_{\{P\in \partial E_3, x_2= {\sqrt{3}/4}\}}  \min( v_1(P), v_2(P))$. 

      We observe now that since $V \leq 0$, $u + 2 C_1\tilde \rho_o^q$ satisfies 
      $$F[u + 2 C_1\tilde \rho_o^q] + h(x)\cdot \nabla (u + 2 C_1\tilde \rho_o^q) |\nabla (u + 2 C_1\tilde \rho_o^q)|^\alpha + V(x) (u + 2 C_1\tilde \rho_o^q)^{1+\alpha } \leq f.$$
       
      The rest of the proof is the same.

  \medskip
\noindent{\em Proof of Corollary \ref{corhold}. }
Suppose that $u$ is a solution in $\Omega $ which contains $B(0, R_o)$.  Let $v$ be defined as $v(x) =: u(Rx)$ . Then $v$ satisfies in $B(0,{R_o\over R})$, 
 $$F(x,\nabla v, D^2v)(x) + Rh(Rx)\cdot \nabla v |\nabla v|^\alpha + R^{2+\alpha}  V(Rx) v^{1+\alpha} = R^{2+\alpha} f(Rx)$$
 Applying Harnack's  inequality for $v$  we get the desired result for $u$.

Let $R_o>0$ such that $B(x_o,4R_o)\subset\Omega^\prime  \subset\subset  \Omega$. We define  for any $R  < R_o$
$$M_i=\max_{B_(x_o,iR)}u,\quad m_i=\min_{B_(x_o,iR)}u$$
for $i=1$ and $i=4$. Then $u-m_i$ is a solution of
$$F[u-m_i]+h(x)\grad (u-m_i)|\grad(u-m_i)|^\alpha=-V(x)u^{1+\alpha}$$
in $B(x_o,iR)$ and hence $u$ satisfies 
$$\sup_{B(x_o,R)} (u(x)-m_4)\leq K\inf_{B(x_o,R)} (u(x)-m_4)+K R^{2+\alpha\over \alpha+1}M_4 |V|_\infty ^{1\over 1+\alpha} $$
In the same way, using the operator $G(x,p,M)=-F(x,p,-M)$, and the function $M_i-u$, we get
\begin{eqnarray*}
G(x,\nabla u, D^2 ( M_i-u)) + h(x)\cdot |\nabla (M_i-u)|^\alpha \nabla (M_i-u)  
&= &V(x) u^{1+\alpha} \
\end{eqnarray*}
in $B(0, iR)$. 
We get  with some constant $K$ which can be taken equal to the previous one   
$$\sup_{B(x_o,R)} (M_4-u(x))\leq K\inf_{B(x_o,R)} (M_4-u(x))+ K R^{2+\alpha\over \alpha+1}M_4 |V|_\infty ^{1\over 1+\alpha} .$$
Summing the inequalities we obtain for some constant $K^\prime $ independant of $R \leq R_o$
$$M_1-m_1\leq\frac{K-1}{K+1} (M_4-m_4)+  K^\prime R^{2+\alpha\over \alpha+1}.  $$
The rest of the proof is classical, just apply Lemma 8.23  in \cite{GT}.

{\em Proof of Corollary \ref{Liou}.}
Let $c_0=\inf_{\R^2} u$ and let $w=u-c_0$. Clearly $w$ satisfies in $\R^2$:
$$
F[w]=0,\ w\geq 0,\ \inf w=0.$$
Suppose by contradiction that $w>0$ somewhere, then applying the strong maximum principle one gets that $w>0$ in the whole of $\R^2$.

By definition of the infimum, for any $\varepsilon>0$ there exists $P\in\R^2$ such that $w(P)\leq\varepsilon$.
Now for any $Q\in\R^2$ consider the ball centered at $P$ and of radius $4|PQ|$, 
by Harnack's inequality and more precisely using Claim 4 in the proof, we get that
$$w(Q)\leq K_2w(P)\leq K_2\varepsilon.$$
Observe that $K_2$ doesn't depend on the distance $|PQ|$  because $h=V\equiv 0$, hence it doesn't depend on the choice of $Q$. 
Since this holds for any $\varepsilon$ we get $w\equiv 0$.

\end{document}